\documentclass[12pt]{amsart}

\usepackage{amsmath}
\usepackage{xcolor}
\usepackage{amsxtra}
\usepackage{amscd}
\usepackage{amsthm}
\usepackage{amsfonts}
\usepackage{amssymb}
\usepackage{eucal}
\usepackage[matrix,arrow,curve]{xy}
\usepackage[dvips]{graphicx}
\usepackage[dvips]{graphics}
\usepackage[T2A]{fontenc}
\usepackage[cp866]{inputenc}
\usepackage{mathtools}
\usepackage{mathrsfs}

\usepackage{hyperref}

\newcommand{\End}{\mathop{\sf End}\nolimits}

\newcommand{\Hom}{\mathop{\sf Hom}\nolimits}

\def \id {{\rm id}}

\theoremstyle{plain}
\newtheorem{Thm}[subsection]{Theorem}
\newtheorem{Cor}[subsection]{Corollary}
\newtheorem{Lem}[subsection]{Lemma}
\newtheorem{Prop}[subsection]{Proposition}
\newtheorem{Conj}[subsection]{Conjecture}
\newtheorem{Ex}[subsection]{Example}

\theoremstyle{definition}
\newtheorem{Def}[subsection]{Definition}

\theoremstyle{remark}

\newtheorem{Rem}[subsection]{Remark}

\newtheorem{conj}{Conjecture}

\numberwithin{equation}{section}

\newif\ifShowLabels
\ShowLabelstrue
\newdimen\theight
\def\TeXref#1{%
    \leavevmode\vadjust{\setbox0=\hbox{{\tt
        \quad\quad  {\small \rm #1}}}%
    \theight=\ht0
    \advance\theight by \lineskip
    \kern -\theight \vbox to
    \theight{\rightline{\rlap{\box0}}%
    \vss}%
    }}%

\ShowLabelsfalse

\renewcommand{\sec}[2]{\section{#2}\label{S:#1}%
    \ifShowLabels \TeXref{{S:#1}} \fi}
\newcommand{\ssec}[2]{\subsection{#2}\label{SS:#1}%
    \ifShowLabels \TeXref{{SS:#1}} \fi}

\newcommand{\sssec}[2]{\subsubsection{#2}\label{SSS:#1}%
    \ifShowLabels \TeXref{{SSS:#1}} \fi}

\newcommand{\refs}[1]{Section ~\ref{S:#1}}
\newcommand{\refss}[1]{Section ~\ref{SS:#1}}
\newcommand{\refsss}[1]{Section ~\ref{SSS:#1}}
\newcommand{\reft}[1]{Theorem ~\ref{T:#1}}
\newcommand{\refl}[1]{Lemma ~\ref{L:#1}}
\newcommand{\refp}[1]{Proposition ~\ref{P:#1}}

\newcommand{\refd}[1]{Definition ~\ref{D:#1}}
\newcommand{\refr}[1]{Remark ~\ref{R:#1}}
\newcommand{\refe}[1]{\eqref{E:#1}}
\newcommand{\refco}[1]{Conjecture ~\ref{Co:#1}}

\newenvironment{thm}[1]%
    { \begin{Thm} \label{T:#1}  \ifShowLabels \TeXref{T:#1} \fi }%
    { \end{Thm} }

\renewcommand{\th}[1]{\begin{thm}{#1} \sl }
\renewcommand{\eth}{\end{thm} }

\newenvironment{lemma}[1]%
    { \begin{Lem} \label{L:#1}  \ifShowLabels \TeXref{L:#1} \fi }%
    { \end{Lem} }
\newcommand{\lem}[1]{\begin{lemma}{#1} \sl}
\newcommand{\elem}{\end{lemma}}

\newenvironment{propos}[1]%
    { \begin{Prop} \label{P:#1}  \ifShowLabels \TeXref{P:#1} \fi }%
    { \end{Prop} }
\newcommand{\prop}[1]{\begin{propos}{#1}\sl }
\newcommand{\eprop}{\end{propos}}

\newenvironment{corol}[1]%
    { \begin{Cor} \label{C:#1}  \ifShowLabels \TeXref{C:#1} \fi }%
    { \end{Cor} }
\newcommand{\cor}[1]{\begin{corol}{#1} \sl }
\newcommand{\ecor}{\end{corol}}

\newenvironment{defeni}[1]%
    { \begin{Def} \label{D:#1}  \ifShowLabels \TeXref{D:#1} \fi }%
    { \end{Def} }
\newcommand{\defe}[1]{\begin{defeni}{#1} \sl }
\newcommand{\edefe}{\end{defeni}}

\newenvironment{remark}[1]%
    { \begin{Rem} \label{R:#1}  \ifShowLabels \TeXref{R:#1} \fi }%
    { \end{Rem} }
\newcommand{\rem}[1]{\begin{remark}{#1}}
\newcommand{\erem}{\end{remark}}

\newenvironment{conjec}[1]%
    { \begin{Conj} \label{Co:#1}  \ifShowLabels \TeXref{Co:#1} \fi }%
    { \end{Conj} }
\renewcommand{\conj}[1]{\begin{conjec}{#1} \sl }
\newcommand{\econj}{\end{conjec}}

\newenvironment{example}[1]%
    { \begin{Ex} \label{Exx:#1}  \ifShowLabels \TeXref{Exx:#1} \fi }%
    { \end{Ex} }
\newcommand{\ex}[1]{\begin{example}{#1} \sl }
\newcommand{\eex}{\end{example}}

\newcommand{\eq}[1]%
    { \ifShowLabels \TeXref{E:#1} \fi
       \begin{equation} \label{E:#1} }
\newcommand{\eeq}{ \end{equation} }

\newcommand{\prf}{ \begin{proof} }
\newcommand{\epr}{ \end{proof} }

\newcommand\alp{\alpha}     

     \newcommand\Gam{\Gamma}
     \newcommand\Del{\Delta}

\newcommand\lam{\lambda}        \newcommand\Lam{\Lambda}

\newcommand\calW{{\mathcal{W}}}

\newcommand\bfk{{\mathbf k}}

\newcommand\LL{\mathbb{L}}
\newcommand\ZZ{\mathbb{Z}}

\newcommand\CC{\mathbb{C}}

 \newcommand\grg{{\mathfrak{g}}}

\newcommand\sdp{\times \hskip -0.3em {\raise 0.3ex
\hbox{$\scriptscriptstyle |$}}} 

\newcommand\Gr{\operatorname{Gr}}

\newcommand\Int{\operatorname{Int}}

\newcommand\x{\times}
\newcommand\ten{\otimes}

\newcommand\nc{\newcommand}

\nc\aff{\operatorname{aff}}
\nc\oGr{\overline{\Gr}}
\nc\Bun{\operatorname{Bun}}
\nc\hgrg{\widehat{\grg}}
\renewcommand\Int{\operatorname{Int}}
\nc\bInt{\overline{\Int}}
\nc\hatLam{\widehat{\Lam}}
\nc\bmu{\overline{\mu}}
\nc\bnu{\overline{\nu}}
\nc\blambda{\overline{\lam}}

\nc\ocalW{\overline{\calW}}
\nc\pos{\operatorname{pos}}
\nc\IH{\operatorname{IH}}
\nc\Rep{\operatorname{Rep}}
\nc\Gal{\operatorname{Gal}}
\nc{\tilGr}{\widetilde{\Gr}}

\nc\Pic{\operatorname{Pic}}
\nc\Prym{\operatorname{Prym}}
\nc\pa{\partial}
\nc\Na{\nabla}

\nc{\HC}{{\mathcal{HC}}}
\nc{\on}{\operatorname}
\nc{\BA}{{\mathbb{A}}}
\nc{\BC}{{\mathbb{C}}}
\nc{\BG}{{\mathbb{G}}}
\nc{\BM}{{\mathbb{M}}}
\nc{\BN}{{\mathbb{N}}}
\nc{\BQ}{{\mathbb{Q}}}
\nc{\BP}{{\mathbb{P}}}
\nc{\BR}{{\mathbb{R}}}
\nc{\BZ}{{\mathbb{Z}}}
\nc{\BS}{{\mathbb{S}}}

\nc{\CA}{{\mathcal{A}}}
\nc{\CB}{{\mathcal{B}}}
\nc{\CalC}{{\mathcal C}}
\nc{\CalD}{{\mathcal D}}
\nc{\CE}{{\mathcal{E}}}
\nc{\CF}{{\mathcal{F}}}
\nc{\CG}{{\mathcal{G}}}
\nc{\CH}{{\mathcal{H}}}
\nc{\CK}{{\mathcal{K}}}
\nc{\CL}{{\mathcal{L}}}
\nc{\CM}{{\mathcal{M}}}
\nc{\CMM}{{\mathcal{M}^{\operatorname{gen}}_\hbar(-\rho)}}
\nc{\CN}{{\mathcal{N}}}
\nc{\CO}{{\mathcal{O}}}
\nc{\CP}{{\mathcal{P}}}
\nc{\CQ}{{\mathcal{Q}}}
\nc{\CR}{{\mathcal{R}}}
\nc{\CS}{{\mathcal{S}}}
\nc{\CT}{{\mathcal{T}}}
\nc{\CU}{{\mathcal{U}}}
\nc{\CV}{{\mathcal{V}}}
\nc{\CW}{{\mathcal{W}}}
\nc{\CX}{{\mathcal{X}}}
\nc{\CY}{{\mathcal{Y}}}
\nc{\CZ}{{\mathcal{Z}}}

\nc{\gen}{{\operatorname{gen}}}
\nc{\cM}{{\check{\mathcal M}}{}}
\nc{\csM}{{\check{\mathcal A}}{}}
\nc{\obM}{{\overset{\circ}{\mathbf M}}{}}
\nc{\oCA}{{\overset{\circ}{\mathcal A}}{}}
\nc{\obA}{{\overset{\circ}{\mathbf A}}{}}
\nc{\ooM}{{\overset{\circ}{M}}{}}
\nc{\osM}{{\overset{\circ}{\mathsf M}}{}}
\nc{\vM}{{\overset{\bullet}{\mathcal M}}{}}
\nc{\nM}{{\underset{\bullet}{\mathcal M}}{}}
\nc{\obD}{{\overset{\circ}{\mathbf D}}{}}
\nc{\cp}{{\overset{\circ}{\mathbf p}}{}}
\nc{\ofZ}{{\overset{\circ}{\mathfrak Z}}{}}

\nc{\fa}{{\mathfrak{a}}}
\nc{\fb}{{\mathfrak{b}}}
\nc{\fg}{{\mathfrak{g}}}
\nc{\fgl}{{\mathfrak{gl}}}
\nc{\fh}{{\mathfrak{h}}}
\nc{\fj}{{\mathfrak{j}}}
\nc{\fm}{{\mathfrak{m}}}
\nc{\fn}{{\mathfrak{n}}}
\nc{\fu}{{\mathfrak{u}}}
\nc{\fp}{{\mathfrak{p}}}
\nc{\frr}{{\mathfrak{r}}}
\nc{\fs}{{\mathfrak{s}}}
\nc{\ft}{{\mathfrak{t}}}
\nc{\fT}{{\mathfrak{T}}}
\nc{\ofT}{{\overline{\mathfrak T}}}
\nc{\ofS}{{\overline{\mathfrak S}}}
\nc{\fsl}{{\mathfrak{sl}}}
\nc{\hsl}{{\widehat{\mathfrak{sl}}}}
\nc{\hgl}{{\widehat{\mathfrak{gl}}}}
\nc{\hg}{{\widehat{\mathfrak{g}}}}
\nc{\chg}{{\widehat{\mathfrak{g}}}{}^\vee}
\nc{\hn}{{\widehat{\mathfrak{n}}}}
\nc{\chn}{{\widehat{\mathfrak{n}}}{}^\vee}

\nc{\fA}{{\mathfrak{A}}}
\nc{\fB}{{\mathfrak{B}}}
\nc{\fD}{{\mathfrak{D}}}
\nc{\fE}{{\mathfrak{E}}}
\nc{\fF}{{\mathfrak{F}}}
\nc{\fG}{{\mathfrak{G}}}
\nc{\fI}{{\mathfrak{I}}}
\nc{\fJ}{{\mathfrak{J}}}
\nc{\fK}{{\mathfrak{K}}}
\nc{\fL}{{\mathfrak{L}}}
\nc{\fM}{{\mathfrak{M}}}
\nc{\fN}{{\mathfrak{N}}}
\nc{\frP}{{\mathfrak{P}}}
\nc{\fS}{{\mathfrak S}}
\nc{\fU}{{\mathfrak{U}}}
\nc{\fZ}{{\mathfrak{Z}}}

\nc{\bb}{{\mathbf{b}}}
\nc{\bc}{{\mathbf{c}}}
\nc{\be}{{\mathbf{e}}}
\nc{\bj}{{\mathbf{j}}}
\nc{\bn}{{\mathbf{n}}}
\nc{\bp}{{\mathbf{p}}}
\nc{\bq}{{\mathbf{q}}}
\nc{\bv}{{\mathbf{v}}}
\nc{\bx}{{\mathbf{x}}}
\nc{\by}{{\mathbf{y}}}
\nc{\bw}{{\mathbf{w}}}
\nc{\bA}{{\mathbf{A}}}
\nc{\bB}{{\mathbf{B}}}
\nc{\bC}{{\mathbf{C}}}
\nc{\bK}{{\mathbf{K}}}
\nc{\bL}{{\mathbf{L}}}
\nc{\bD}{{\mathbf{D}}}
\nc{\bH}{{\mathbf{H}}}
\nc{\bM}{{\mathbf{M}}}
\nc{\bN}{{\mathbf{N}}}
\nc{\bS}{{\mathbf{S}}}
\nc{\bT}{{\mathbf{T}}}
\nc{\bV}{{\mathbf{V}}}
\nc{\bW}{{\mathbf{W}}}
\nc{\bX}{{\mathbf{X}}}
\nc{\bP}{{\mathbf{P}}}
\nc{\bZ}{{\mathbf{Z}}}

\nc{\sA}{{\mathsf{A}}}
\nc{\sB}{{\mathsf{B}}}
\nc{\sC}{{\mathsf{C}}}
\nc{\sD}{{\mathsf{D}}}
\nc{\sF}{{\mathsf{F}}}
\nc{\sK}{{\mathsf{K}}}
\nc{\sM}{{\mathsf{M}}}
\nc{\sO}{{\mathsf{O}}}
\nc{\sQ}{{\mathsf{Q}}}
\nc{\sP}{{\mathsf{P}}}
\nc{\sV}{{\mathsf{V}}}
\nc{\sW}{{\mathsf{W}}}
\nc{\sZ}{{\mathsf{Z}}}
\nc{\sfp}{{\mathsf{p}}}
\nc{\sr}{{\mathsf{r}}}
\nc{\sfb}{{\mathsf{b}}}
\nc{\sfc}{{\mathsf{c}}}
\nc{\sd}{{\mathsf{d}}}
\nc{\sg}{{\mathsf{g}}}
\nc{\sfl}{{\mathsf{l}}}

\nc{\BK}{{\bar{K}}}

\nc{\tA}{{\widetilde{\mathbf{A}}}}
\nc{\tB}{{\widetilde{\mathcal{B}}}}
\nc{\tg}{{\widetilde{\mathfrak{g}}}}
\nc{\tG}{{\widetilde{G}}}
\nc{\TM}{{\widetilde{\mathbb{M}}}{}}
\nc{\tO}{{\widetilde{\mathsf{O}}}{}}
\nc{\tU}{{\widetilde{\mathfrak{U}}}{}}
\nc{\TZ}{{\tilde{Z}}}
\nc{\tZ}{\widetilde{Z}{}}
\nc{\tx}{{\tilde{x}}}
\nc{\tbv}{{\tilde{\bv}}}
\nc{\tfP}{{\widetilde{\mathfrak{P}}}{}}
\nc{\tz}{{\tilde{\zeta}}}
\nc{\tmu}{{\tilde{\mu}}}

\nc{\td}{\ddot{\underline{d}}{}}
\nc{\tzeta}{\widetilde{\zeta}{}}
\nc{\hd}{{\widehat{\underline{d}}}}
\nc{\hG}{{\widehat{G}}}
\nc{\hBP}{\widehat{\mathbb P}{}}
\nc{\hQ}{{\widehat{Q}}}
\nc{\hsM}{\widehat{\mathsf M}{}}
\nc{\hfM}{\widehat{\mathfrak M}{}}
\nc{\hCP}{\widehat{\mathcal P}{}}
\nc{\hCR}{\widehat{\mathcal R}{}}
\nc{\hCS}{{\widehat{\mathcal S}}}
\nc{\hfZ}{\widehat{\mathfrak Z}{}}

\nc{\urho}{\underline{\rho}}
\nc{\uB}{\underline{B}}
\nc{\uC}{{\underline{\mathbb{C}}}}
\nc{\uk}{{\underline{\bfk}}}
\nc{\ui}{\underline{i}}
\nc{\ofP}{{\overline{\mathfrak{P}}}}

\nc{\hrho}{{\hat{\rho}}}

\nc{\unl}{\underline}
\nc{\ol}{\overline}
\nc{\one}{{\mathbf{1}}}
\nc{\two}{{\mathbf{t}}}

\nc{\Tot}{{\mathop{\operatorname{\rm Tot}}}}
\nc{\Hilb}{{\mathop{\operatorname{\rm Hilb}}}}
\nc{\CHom}{{\mathop{\operatorname{{\mathcal{H}}\it om}}}}
\nc{\defi}{{\mathop{\operatorname{\rm def}}}}
\nc{\length}{{\mathop{\operatorname{\rm length}}}}

\nc{\Cliff}{{\mathsf{Cliff}}}
\nc{\Fl}{{\mathsf{Fl}}}
\nc{\Fib}{{\mathsf{Fib}}}
\nc{\Coh}{{\mathsf{Coh}}}
\nc{\FCoh}{{\mathsf{FCoh}}}

\nc{\reg}{{\text{\rm reg}}}

\nc{\cplus}{{\mathbf{C}_+}}
\nc{\cminus}{{\mathbf{C}_-}}
\nc{\cthree}{{\mathbf{C}_*}}
\nc{\Qbar}{{\bar{Q}}}
    
\nc{\bh}{{\bar{h}}}
\nc{\bOmega}{{\overline{\Omega}}}
\nc\tGr{\widetilde{\Gr}}

\nc{\seq}[1]{\stackrel{#1}{\sim}}
\nc\ogu{\overline{G/U}}
\nc\chlam{\check{\lam}}

\nc\St{\operatorname{St}}

\nc\uS{\underline{S}}
\nc\QM{\mathcal{QM}}
\nc\FT{\mathsf{FT}}
\nc{\Ca}{\underline{C_a}}
\nc{\SCa}{\underline{S^mC_a}}
\nc{\sCa}{\Ca\x_\CA\underline{S^{m-1}C_a}}

\nc\Shv{\operatorname{Shv}}
\nc\cmod{\textup{-comod}}
\nc\lmod{\textup{-mod}}
\nc\rmod{\textup{mod-}}
\nc\bimod{\textup{-mod-}}
\nc\mon{\operatorname{mon}}
\nc\pt{\textup{pt}}
\nc\nil{\textup{nil}}
\nc\RG{\textup{R}\Gamma}
\nc\RGc{\textup{R}\Gamma_{c}}
\nc\fmo{free-monodromic }

\nc\fra{\mathfrak{a}}
\nc\frg{\mathfrak{g}}
\nc\frb{\mathfrak{b}}
\nc\frn{\mathfrak{n}}
\nc\frh{\mathfrak{h}}
\nc\frp{\mathfrak{p}}
\nc\frk{\mathfrak{k}}

\nc\Lot{\stackrel{\LL}{\otimes}}

\newcommand\quash[1]{}
\renewcommand\c\circ

\renewcommand\a\alpha
\renewcommand\b\beta
\newcommand\g\gamma
\newcommand\G\Gamma
\renewcommand\d\delta
\newcommand\D\Delta

\nc\io{\iota}

\renewcommand\r\rho

\nc\coker{\textup{coker}}
\nc\ev{\textup{ev}}
\nc\pro{\textup{pro}}
\nc\perf{\textup{perf}}
\nc\Perf{\textup{Perf}}
\nc\Tilt{\textup{Tilt}}
\nc\triv{\textup{triv}}
\nc\Lie{\textup{Lie}}
\nc\bt{\boxtimes}
\nc\SU{\textup{SU}}
\nc\PU{\textup{PU}}
\nc\SBim{\textup{SBim}}
\nc\TGR{\Tilt(\CM_{G_{\BR}})}
\nc\THG{\Tilt(\CH_{G})}
\nc\Frac{\textup{Frac}}
\nc\red{\textup{red}}
\nc\Spf{\textup{Spf}\ }
\nc\Inv{\textup{Inv}}
\nc\Stab{\textup{Stab}}

\nc\bstar{\overline{\star}}
\nc\LS{\textup{LS}}
\nc\act{\textup{act}}

\title{McKay correspondence and orbifold equivalence}
\dedicatory{In memory of S.\,Natanzon}
\author{Andrei Ionov}

\address{
Boston College, Department of Mathematics, Maloney Hall, Fifth Floor,
Chestnut Hill, MA
02467-3806, United States }
\email{ionov@bc.edu}

\begin{document}

\maketitle

\begin{abstract}
We prove that a pair of singularities related by a transformation arising from the McKay correspondence are orbifold equivalent. From this we deduce a new proof of a McKay type equivalence for the matrix factorization categories.
\end{abstract}

\sec{}{Introduction}

For an algebraic variety $V$ with an action of a finite group $G$ and a crepant resolution $Y$ of the quotient $V/G$ the {\itshape McKay correspondence} refers to various statements relating the $G$-equivariant geometry of $V$ and the geometry of $Y$. The categorified version provides an equivalence of categories $$D^G(\Coh(V))\cong D(\Coh(Y))$$ under various assumptions on $V, G$ and $Y$. For example, in case of $\dim V=2,3$ and $Y$ being a $G$-Hilbert scheme of $V$ the result of this type were established in \cite{BGK} and in \cite{BK} this was proved  for the case of symplectic resolutions.

Motivated by conformal field theories and defined in \cite{CR2}, the {\itshape orbifold equivalence} is an equivalence relation on the set of polynomial potentials $W\colon \BC^n\to\BC$. Namely, $W_1$ and $W_2$ are equivalent if there exists a matrix factorization of $W_2-W_1$, whose {\itshape left} and {\itshape right quantum dimensions} in the sense of bicategories are invertible (see also \cite{CR}, \cite{CM}, \cite{CRCR}). 

To unify these two concepts we associate to a $G$-invariant function $f$ on $V$ a function $\hat{f}$ on $Y$. In the present paper we prove:
\th{}(\reft{main})
Under the assumptions of \refs{set} the singularities of $f$ and $\hat{f}$ are orbifold equivalent. 
\eth
The observation of this sort were made in \cite{BT}, which was a motivation for the present paper. It was noticed there that several of the known examples of orbifold equivalence fit into this setting. Note that in the most of the known cases the orbifold equivalences are constructed by the computational methods (\cite{CRCR}, \cite{NRC}, \cite{RW}, \cite{KRC}). Our methods are geometric and allow to reprove some of the known results. Studying the examples makes us hope that a class of orbifold equivalences given by a matrix factorization $X$, such that the product of the quantum dimensions is in $\BQ^{>0}$, is generated, in a sense, by the equivalences provided by the McKay correspondence (\refco{conj}). 

From the orbifold equivalence we deduce a new proof of the McKay type equivalence of matrix factorizations categories proved in the most general setting in \cite{BP} (\reft{mod}).

Our assumptions mainly concern the McKay correspondence itself. They are satisfied in the setting of \cite{BGK} and it was communicated to the author by R.\,Bezrukavnikov that they are expected to hold in the setting of \cite{BK}, although they are not known to hold in general. In particular, in the essential special case of symplectic resolutions, the Hilbert scheme of points on $\BC^2$ resolving the singularity of $(\BC^2)^n/S_n$, they are known to hold (see \cite{H}).

It would be interesting to see what do our results mean in the language of the conformal field theories. In particular, it seems natural to ask if there is a "physical" or "categorical" version of the McKay correspondence, which could provide orbifold equivalences in a wider range of bicategories. We refer to \cite{R} for the discussion of Landau-Ginzburg/conformal field theory correspondence, which should be a key to these questions.

The paper is structured as follows. In \refs{back} we briefly recall the background on matrix factorizations and orbifold equivalence. In \refs{set} we introduce our setting, including the discussion of the McKay correspondence, and state our results in details. In \refs{proof} we prove our results. In \refs{exmpl} we will look at the known examples of the orbifold equivalences and see how they fit into our setting and conjecture what type of orbifold equivalences are generated by McKay correspondence.

\ssec{}{Conventions and notations} We will work over the field of complex numbers $\CC$. However, nothing changes if one replaces $\CC$ by any algebraically closed field of characteristic $0$.

For an algebraic variety $S$ we put $D(S)$ for the bounded derived category of coherent sheaves on $S$ and $\Perf(S)$ for the category of perfect complexes. If there is an action of a group $G$ on $S$ we use $D^G(S)$ for the bounded derived category of equivariant sheaves. All functors are assumed to be derived unless the opposite is specified. 

We put $K_S$ for the canonical line bundle of $S$ and if $\phi\colon S_1\to S_2$ we put $K_\phi=K_{S_1/S_2}$ for the relative canonical sheaf. 

\ssec{}{Acknowledgements} The author is grateful to A.\,Basalev, R.\,Bezrukavnikov, A.\,Kuznetsov, S.\,Natanzon, A.\,Petrov, A.\,Ros Camacho and M.\,Wemyss for useful discussions. The author is grateful to A.\,Basalev, F.\,Kogan and A.\,Utiralova for reading the original manuscript of the paper. The author wants to express gratitude  to the anonymous referee for helpful comments and suggestions. The author was funded by RFBR, project number 19-31-90078.


\sec{back}{Background}

\ssec{}{Matrix factorizations} We refer to \cite{O1}, \cite{O2}, \cite{D} and \cite{CM} for a comprehensive discussion of the subject in different contexts. 

Let $S$ be an algebraic variety and let $W\in\Gam(\CO_S)$ be a global function on $S$. Respectively, let $R$ be a commutative ring and $W\in R$ be its element.

\defe{}
The category of matrix factorizations $\mathrm{MF}(W)=\mathrm{MF}(S,W)$ (respectively $\mathrm{MF}(R,W)$) of $W$ is a differential $\BZ/2\BZ$-graded category, whose objects are pairs $(X,d_X)$, where $X=X^0\oplus X^1$ is a locally free $\BZ/2\BZ$-graded sheaf (of finite rank) on $S$ (respectively projective module over $R$) and $d_E$ is its odd endomorphism satisfying $d_X^2=W\cdot \id_X$. 

The morphism complex between $(X,d_X)$ and $(X',d_{X'})$ is a graded space $\Hom(X,X')$ of all maps between $X$ and $X'$ with a differential given by $$\partial(\alp)=d_{X'}\a-(-1)^{|\a|}\a d_X.$$ 

We put $\mathrm{HMF}(W)=\mathrm{HMF}(S,W)$ (respectively $\mathrm{HMF}(R,W)$) for the idempotent completion of the corresponding homotopy category. 

\edefe

One easily adds additional structures, e.g. grading or equivariance, into the definition. 

Let $$D_{sing}(\{W=0\}):=D(\{W=0\})/\Perf(\{W=0\})$$ be the singularity category of $W$. The following is the key result on the matrix factorizations:

\th{O}(\cite[Theorem 3.5]{O2}) There is an equivalence of triangulated categories:

\eq{O} D_{sing}(\{W=0\})\cong\mathrm{HMF}(W).\eeq

\eth


Let $S_1, S_2, S_3$ be three varieties together with three functions $W_1, W_2, W_3$ on them. 
Then the tensor product over $\CO_{S_2}$ induces a functor \eq{ten}-\ten_{\CO_{S_2}}-\colon\mathrm{HMF}(S_3\x S_2,W_3-W_2)\x\mathrm{HMF}(S_2\x S_1,W_2-W_1)\to\eeq $$\to\mathrm{HMF}(S_3\x S_1,W_3-W_1).$$ Note that there is a technical difficulty that could be overcome: naively, the tensor product will be a matrix factorization of infinite rank. This issue could be resolved using, for example \cite[Section 12]{DM}. 

Note also that we can plug $S_i=pt, W_i=0$ for some $i$ into \refe{ten}. For $S_1=pt, S_2=S_3=S$ and $W_1=0, W_2=W_3=W$ put \eq{double}\widetilde{W}:=W_3-W_2=W\ten1-1\ten W\in\Gam(\CO_S)\otimes_\BC\Gam(\CO_S)=\Gam(\CO_{S\x S}).\eeq Then there is the unique up to homotopy equivalence object $\Del_W$ of $\mathrm{HMF}(S\x S,\widetilde{W})$ under \refe{ten} representing the identity functor on $\mathrm{HMF}(S,W)$ and playing the role of the monoidal unit. 

The tensor product and objects $\Del_W$ satisfy all the naturally expected properties.

\ssec{orb}{Orbifold equivalence} We refer to \cite{CR}, \cite{CM}, \cite{CR2} and \cite{CRCR} for the detailed discussion of this topic.


Let $R_1, R_2$ be either two polynomial rings or two rings of formal power series, such that the number of variables has the same parity. In the latter case one should work with the completed tensor products. The localization result \cite[Theorems 4.11, 5.7]{D} allows one to switch between these settings and transfer the results. 

Pick $W_1\in R_1, W_2\in R_2$ and let $X=(X,d_X)$ be in $\mathrm{HMF}(R_2\ten_\BC R_1, W_2-W_1)$. We define $X^\vee=(X^\vee, d_{X^\vee})$ to be the dual matrix factorization given by $X^\vee=\Hom_{R_2\ten_\BC R_1}(X,R_2\ten_\BC R_1)$ and $d_{X^\vee}\alp=(-1)^{1+|\a|} \a\circ d_X$. This is an object of $\mathrm{HMF}(R_1\ten_\BC R_2, W_1-W_2)$. Shifting it by the parity of the number of variables we obtain an object $X^\dag$ of $\mathrm{HMF}(R_1\ten_\BC R_2, W_1-W_2)$. The functors $X\ten_{R_1}-$ and $X^\dag\ten_{R_2}-$ are two sided adjoint (\cite[Section 7]{CM}). Note that it is here that we use the fact that the numbers of variables in the two considered rings have the same parity, as otherwise the shift by $1$ will appear in the adjunction. The evaluation and coevaluation maps are induced by $4$ morphisms in $\mathrm{HMF}(\widetilde{W_1})$ and $\mathrm{HMF}(\widetilde{W_2})$:
$$\Del_{W_1}\to X^\dag\ten_{R_2}X\to \Del_{W_1},$$
$$\Del_{W_2}\to X\ten_{R_1} X^\dag\to \Del_{W_2}.$$
Taking the compositions one obtains elements $$\dim_l(X)\in\End(\Del_{W_1})=\mathrm{Jac}({W_1})$$ and $$\dim_r(X)\in\End(\Del_{W_2})=\mathrm{Jac}({W_2}),$$ where $\mathrm{Jac}({W})$ is the Jacobi ring of $W$. The elements $\dim_l(X)$ and $\dim_r (X)$ are respectively called {\itshape left} and {\itshape right quantum dimensions}.
 
By \cite[Proposition 8.5]{CM} the quantum dimensions enjoy the following properties:
\eq{1}
\dim_l(\Del_W)=\dim_r(\Del_W)=1;
\eeq
\eq{2}
\dim_l(X)=\dim_r(X^\dag);
\eeq
\eq{3}
\dim_l(X_1\ten X_2)=\dim_l(X_2)\dim_l(X_1), \enskip \dim_r(X_1\ten X_2)=\dim_r(X_1)\dim_r(X_2).
\eeq
The quantum dimensions can be computed by the direct formulas in terms of $d_X$ given by \cite[Proposition 8.2]{CM}. 

\defe{orb}
We say that $W_1\in R_1$ and $W_2\in R_2$ are orbifold equivalent and write $W_1\sim W_2$ if there is an object $X$ of $\mathrm{HMF}(R_2\ten_\BC R_1, W_2-W_1)$ with invertible left and right quantum dimensions.
\edefe

Given a grading on rings $R_1$ and $R_2$, such that $W_1$ and $W_2$ are of pure degree, in \refd{orb} one may additionally ask for $X$ to be graded, which provides the graded version of the orbifold equivalence. In the case of $\BQ^{\ge 0}$-grading, such that the $0$ degree component of $R_1$ and $R_2$ are equal to $\BC\cdot1$, the quantum dimension is known to be an element of $\BC$ (\cite[Lemma 2.5]{CRCR}). This localization procedure can always be done. 


The orbifold equivalence provides an important relation between the categories $\mathrm{HMF}(W_1)$ and $\mathrm{HMF}(W_2)$. The endofunctor 
$$X^\dag\ten_{R_2}X\ten_{R_1}-\colon \mathrm{HMF}(W_1)\to \mathrm{HMF}(W_1)$$
by adjunction admits a structure of a monad, i.e. it has a multiplication structure given by the composition
$$X^\dag\ten_{R_2}X\ten_{R_1}X^\dag\ten_{R_2}X\to X^\dag\ten_{R_2}\Del_{W_2}\ten_{R_2}X\cong X^\dag\ten_{R_2}X.$$ It makes sense to talk about objects in $\mathrm{HMF}(W_1)$ with an action of $X^\dag\ten_{R_2}X$. We denote by $X^\dag\ten_{R_2}X-\bmod_{\mathrm{HMF}(W_1)}$ the category of such module objects. We have

\th{monad}(\cite[Corollary 1.6]{CRCR}, see also \cite{CR2})
There is an equivalence of triangulated categories 
$$\mathrm{HMF}(W_2)\cong X^\dag\ten_{R_2}X-\bmod_{\mathrm{HMF}(W_1)}.$$
\eth

\sssec{var}{Orbifold equivalence for functions on varieties} We need to also make sense of the orbifold equivalence in a more general setting of a function $W$ on an algebraic variety $S$. Assume that $W$ has an isolated singularity at a unique smooth point $s\in S$ and is smooth everywhere else. Assume that there is an affine chart isomorphic to an affine space $U\subset S$ and containing $s$. We can then restrict $W$ to $U$ to get an element $W|_U$ of the polynomial ring $\BC[U]$. We will discuss the existence of such affine charts in the relevant situations later in \refss{ex}. Alternatively, we can consider the completion of the local ring $\widehat{\CO}_{S,s}$ of $S$ at $s$. It is a formal power series ring and the germ $W_s$ of $W$ at $s$ provides its element. This localization procedure can always be done. The localization equivalences of the matrix factorization categories \cite[Proposition 1.14]{O1} and \cite[Theorem 4.11, 5.7]{D}  provide canonical equivalences $$\mathrm{HMF}(S,W)\cong\mathrm{HMF}(\BC[U],W|_U)\cong\mathrm{HMF}(\widehat{\CO}_{S,s},W_s).$$ If $(S_1, W_1)$ and $(S_2, W_2)$ are two such pairs we similarly have 
$$\mathrm{HMF}(S_1\x S_2,W_1+W_2)\cong$$ $$\cong\mathrm{HMF}(\BC[U_1]\ten_\BC\BC[U_2],W_1|_{U_1}+W_2|_{U_2})\cong\mathrm{HMF}(\widehat{\CO}_{S,s_1}\widehat{\ten}_\BC\widehat{\CO}_{S,s_2},W_{s_1}+W_{s_2}),$$
where $s_i$ are the unique singularity points and $U_i$, if they exist, are two suitable affine charts. The equivalences are compatible with the tensor products. This implies that we can canonically fit $\mathrm{HMF}(S,W)$ for $W$ with isolated singularity at the unique point into the bicategory of polynomial potentials (in case of the existence of the affine chart) or power series potentials (in all cases) and their matrix factorizations. In particular, we can talk about orbifold equivalence between such functions. Canonicity of the equivalences implies that the choices made at this point will not impact our results.


\sec{set}{Setting and the formulation of the main result} 


\ssec{}{McKay correspondence} 
Let $V$ be a smooth quasiprojective complex variety with an action of a finite group of automorphisms $G$. Assume further that the canonical bundle $K_V$ of $V$ is locally trivial as a $G$-sheaf, i.e. if for each point $x\in V$ the stabilizer of $x$ in $G$ acts on the tangent space $T_xV$ by an element with determinant $1$. Let $M=V/G$ be the quotient, which then has only Gorenstein singularities and $\pi\colon V\to M$ be the (finite) quotient map.

Assume that there is a resolution of singularties $\tau\colon Y\to M$, which is crepant, i.e. $\tau^*(K_M)=K_Y$. Let $Z=Y\x_M V\subset Y\x V$ be the Cartesian product so that we have the diagram

\eq{cart}
\xymatrix{
 & Z \ar[dl]_p \ar[dr]^q &  \\
 Y  \ar[dr]_\tau & & V\ar[dl]^\pi   \\
 & M &. 
}
\eeq
Note that there is a $G$ action on $Z$ coming from the action on the second factor such that $Y=Z/G$. The quotient map $p$ is finite and we further assume it to be flat \footnote{Note that since $Y$ is smooth it is sufficient to know that $Z$ is Cohen-Macaulay.}. 

Fix a locally free sheaf $\CE$ on $Z$, which is a pull back along $p$ and consider a Fourier-Mukai transform with kernel $\CE$:
\eq{}
\Phi_\CE^G\colon D(Y)\to D^G(V)
\eeq
given by
\eq{}
\Phi_\CE^G(-):=q_*(p^*(-\ten\rho_0)\ten_{\CO_Z}\CE),
\eeq
where $-\ten\rho_0\colon D(Y)\to D(Y)^G$ is the functor of taking the trivial $G$-representation (the $G$-action on $Y$ is trivial). Let $\CF=\CE^\vee$ be the dual sheaf, so that we have a right adjoint Fourier-Mukai transform:
\eq{}
\Phi_{\CE^\vee}^G\colon D^G(V)\to D(Y)
\eeq
given by
\eq{}
\Phi_{\CE^\vee}^G(-):=(p_*(q^!(-)\ten_{\CO_Z}\CE^\vee))^G,
\eeq
where $(-)^G\colon D(Y)^G\to D(Y)$ is the functor of taking $G$-invariants. We will assume that $\Phi^G_\CE$ and $\Phi^G_{\CE^\vee}$ are mutually inverse equivalences.

\ex{} The situation considered in \cite{BGK} fits into the present setting: in this case $\CE=\CO_Z=p^*\CO_Y$ and the other conditions are stated directly in {\itshape loc. cit.}
\eex

Let $\Phi_{\CE}(-)=q_*(p^*(-)\ten_{\CO_Z}\CE)$ and $\Phi_{\CE^\vee}(-):=p_*(q^!(-)\ten_{\CO_Z}\CE^\vee)$ be the nonequivariant version of the adjoint pair above. We will make use of the following observation:

\lem{comp}
We have an equivalence of functors $\Phi_{\CE^\vee}\circ\Phi_{\CE}(-)\cong p_*\CO_Z\ten_{\CO_Y}-$.
\elem

\prf

By \cite[Theorem 1.3]{N} the functor $p^*$ induces an equivalence between $D(Y)$ and the full subcategory $I^G(Z)\subset D^G(Z)$ consisting of sheaves, such that for each $z\in Z$ the stabilzer $G_z\subset G$ acts trivially on the stalk of the sheaf at $z$. The inverse functor is given by the composition of $p_*$ and $(-)^G$. Note that by our assumptions $-\ten_{\CO_Z}\CE$ and $-\ten_{\CO_Z}\CE^\vee$ preserves $I^G(Z)$ and so does $-\ten_{\CO_Z} K_{q}$.

Given $\CF\in D(Y)$ we have $p^*\CF\ten_{\CO_Z}\CE$ is in $I^G(Z)$ and so it is isomorphic to $p^*\CF'$ for $\CF'\cong (p_*p^*\CF\ten_{\CO_Z}\CE)^G$ on $Y$. 
Since $G$ action on $Z=Y\x_M V$ comes from the action on the second factor, $q^*q_*$ preserves $I^G(Z)$. In particular $q^*q_*p^*\CF'$ is in $I^G(Z)$ and so is $$q^!q_*p^*\CF'\ten_{\CO_Z}\CE^\vee\cong q^*q_*p^*\CF'\ten_{\CO_Z} K_{q}\ten_{\CO_Z}\CE^\vee.$$ We conclude that 
$$q^!\Phi_{\CE}(\CF)\ten_{\CO_Z}\CE^\vee \cong p^*(p_*(q^!\Phi_{\CE}(\CF)\ten_{\CO_Z}\CE^\vee))^G\cong p^*\CF.$$ The projection formula implies 
$$\Phi_{\CE^\vee}\circ\Phi_{\CE}(\CF)\cong p_*p^*\CF\cong p_*\CO_Z\ten_{\CO_Y}\CF.$$
Finally, all the equivalences used are functorial.

\epr

\ssec{}{Introducing the singularity into the picture}

Let $f\colon V\to \BC$ be a $G$-invariant function. It then descends to a function on $M$ and by taking composition with $\tau$ we obtain a function $\hat{f}\colon Y\to\BC$. Let $f$ and $\hat{f}$ both have isolated singularities at the unique points $v\in V$ and $y\in Y$ respectively. 

The following statement was proved in the most general setting in \cite[Theorem 1.1]{BP}. In the present setting we will give a new proof in the next section.

\th{mod}
There is an equivalence of triangulated categories
$$\mathrm{HMF}(\hat{f})\cong\mathrm{HMF}^G(f).$$
\eth



We will work in the context of \refsss{var} and fix the respective choices for $(V,f)$ and $(Y,\hat{f})$. Note that $v$ is necessary $G$-fixed and so $G$ induces a linear action on the local ring $\widehat{\CO}_{V,v}$. It, therefore, make sense to talk about $\mathrm{HMF}^G(\widehat{\CO}_{V,v},f_s)$, which is equivalent to $\mathrm{HMF}^G(V,f)$. In case of working with an affine space chart $v\in U\subset V$ we should further assume that it is $G$-invariant. We shall prove

\th{main}
There exists a matrix factorization of $\hat{f}-f$ as a function on $Y\x V$, which after the localization at $(y,v)$ has invertible left and right quantum dimensions.
\eth

\ssec{}{Grading} 

Suppose additionally that there is an action of $\BC^\x$ on $V$ commuting with the action of $G$ and an action $\BC^\x$ on $Y$, altogether making \refe{cart} to be the $\BC^\x$-equivariant Cartesian diagram. In this case we further suppose that $\CE$ admits a structure of $\BC^\x$-equivariant sheaf, which we fix, so that $\Phi_\CE, \Phi_{\CE^\vee}, \Phi^G_\CE, \Phi^G_{\CE^\vee}$ define functors between the $\BC^\x$-equivariant categories.

\ex{}
In the setting of \cite{BGK} the $\BC^\x$-action on $V$ commuting with the $G$-action defines the compatible $\BC^\x$-action on $Y$.  Indeed, recall that $Y$ is an irreducible component of $G$-$\Hilb$. Then if $C\subset V$ is a $G$-cluster and $t\in\BC^\x$, then $t(C)\subset V$ is also a $G$-cluster and as $\BC^\x$ is connected we stay in the same component. The sheaf $\CE=\CO_Z$ is naturally a $\BC^\x$-equivariant sheaf.
\eex

Let $f$ be homogenous with respect to the $\BC^\x$-action, i.e. for a fixed $D\in \ZZ$ and any $t\in\BC^\x$ we have $t^*(f)=t^Df$. The function $\hat{f}$ is then also homogenous with the same value of $D$. This implies, in particular, that $v$ and $y$ are the fixed points of the action and we have a $\ZZ$-grading on the corresponding local rings, such that the images of $f$ and $\hat{f}$ are of pure degree $D$.

The following proposition makes the matrix factorization we construct into the graded one.

\prop{grad}
In the present setting the matrix factorization of \reft{main} could be chosen to be $\BC^\x$-equivariant ($\BZ$-graded).
\eprop

\sec{proof}{Proof of the main result}

Note that by construction $\hat{f}-f$ vanishes on $Z$. Thus, we can view $\CE$ as a coherent sheaf on the zero variety $\{\hat{f}-f=0\}\subset Y\x V$. We can now project $\CE$ to the singularity category $D_{sing}(\{\hat{f}-f=0\})$. Now \reft{O} yields a matrix factorization $X_\CE$ of $\hat{f}-f$ associated to $\CE$. In the same way there is a matrix factorization $X_{\CE^\vee}$ of $f-\hat{f}$, which is two sided adjoint to $X_\CE$. Recall the notation \refe{double}.

\prop{iso}
In the category $\mathrm{HMF}(\widetilde{\hat{f}})$ we have $$X_{\CE^\vee}\ten_{\CO_V} X_\CE\cong\BC[G]\ten_\BC\Del_{\hat{f}},$$ where $\BC[G]$ is the vector space of functions on $G$.
\eprop

\prf[\refp{iso} implies \reft{main} and \refp{grad}]

By
\refe{3} we have $$\dim_l(X_\CE)\dim_l(X_{\CE^\vee})=\dim_l(\BC[G]\ten_\BC\Del_{\hat{f}}),$$ $$\dim_r(X_{\CE^\vee})\dim_r (X_\CE)=\dim_r(\BC[G]\ten_\BC\Del_{\hat{f}}).$$ By \refe{1} 
and the fact that $\dim_\BC\BC[G]=|G|$ we have $$\dim_l(\BC[G]\ten_\BC\Del_{\hat{f}})=\dim_r(\BC[G]\ten_\BC\Del_{\hat{f}})=|G|.$$ Since $|G|\ne 0$, the statement follows.

In the setting of \refp{grad} note that $\{\hat{f}-f=0\}$ is preserved by the $\BC^\x$-action and the equivalence \refe{O} is compatible with the actions, which allows to pass to the $\BC^\x$-equivariant categories. Since $\CE$ was a $\BC^\x$-equivariant sheaf $X_\CE$ is an equivariant matrix factorization.

\epr

\rem{dim}
As $X^{\dag}_\CE\cong X_{\CE^\vee}$ and \refe{1} holds, we have proven that $\dim_l(X_\CE)\dim_r(X_{\CE})=|G|$.
\erem

\prf[Proof of \refp{iso}] By definition, the tensor product functors $X_\CE\ten_{\CO_Y}-$ and $X_{\CE^\vee}\ten_{\CO_V}-$ are  the Fourier-Mukai transforms with the kernels $X_\CE$ and $X_{\CE^\vee}$ respectively. It follows that under the equivalence \refe{O} these functors respectively match the functors induced by $\Phi_\CE$ and $\Phi_{\CE^\vee}$. 

Now \refl{comp} implies that under the equivalence \refe{O} applied to $W=f$ 
the functor $X_{\CE^\vee}\ten_{\CO_V} X_\CE\ten_{\CO_Y}-$ on the right hand side matches the functor induced by $p_*\CO_Z\ten_{\CO_Y}-$ on the left. Since $p$ is flat and finite, the sheaf $p_*\CO_Z$ is locally free. Moreover, its fiber over a point is identified with $\BC[G]$. As we work locally near $y$ the sheaf $p_*\CO_Z$ is isomorphic to the trivial sheaf with the fiber $\BC[G]$. We see that the functor induced by $p_*\CO_Z\ten_{\CO_Y}-$ on the singularity category corresponds to $\BC[G]\ten_\BC\Del_{\hat{f}}\ten_{\CO_Y}-$ on the matrix factorization category. The statement follows.
\epr

\prf[Proof of \reft{mod}]\footnote{The idea of this argument was suggested to the author by the anonymous referee.} We can upgrade the proof and the statement of \refp{iso} to the isomorphism of monads, with the monad structure on $\BC[G]\ten_\BC\Del_{\hat{f}}\ten_{\CO_Y}-$ being the one coming from the algebra structure of $\BC[G]$. Indeed, the argument of \refp{iso} implies that the monad structure of interest is identified with the monad structure on $p_*\CO_Z\ten_{\CO_Y}-$ or in other words $\CO_Y$-algebra structure on $p_*\CO_Z\simeq p_*p^*\CO_Y$, which is fiberwise isomorphic to the algebra $\BC[G]$ of functions on $G$. 

Applying \reft{monad} we obtain

$$\mathrm{HMF}(f)\cong X_{\CE^\vee}\ten_{\CO_V} X_\CE-\bmod_{\mathrm{HMF}(\hat{f})}\cong \BC[G]\ten_\BC\Del_{\hat{f}}-\bmod_{\mathrm{HMF}(\hat{f})}.$$
Since $\Del_{\hat{f}}$ acts on $\mathrm{HMF}(\hat{f})$ as identity, the latter category is just the category of matrix factorizations of $\hat{f}$ with an action of $\BC[G]$ commuting with the differential (equivalently, $G$-graded matrix factorizations), which we denote by just $\BC[G]-\bmod_{\mathrm{HMF}(\hat{f})}$.
As the functors involved were compatible with $G$-action, the equivalence above is the equivalence of $G$-categories, where $G$ acts on the right hand side through its action on $V$ and on the left hand side through the translation action on $\BC[G]$. Passing to the categories of $G$-equivariant objects and using the fact that $\BC[G]^G\cong\BC$ we conclude 
$$\mathrm{HMF}^G(f)\cong (\BC[G]-\bmod_{\mathrm{HMF}(\hat{f})})^G\cong\mathrm{HMF}(\hat{f})$$
as desired.

\epr


\sec{exmpl}{Relation to the known examples and conjecture} 

\ssec{ex}{Examples}

\sssec{}{$A_{2b-1}\sim D_{b+1}$ in \cite{CR2} and \cite{CRCR}} In the notations of {\itshape loc. cit.} we take $G=\BZ/2\BZ$ acting by $(u,v)\mapsto (-u,-v)$ so that $f$ is of type $A_{2b-1}$ and the singularity of $\hat{f}$ is of type $D_{b+1}$. This was observed in \cite{BTW1}, where  it is also shown that one can chose a suitable affine chart on $Y$ isomorphic to the affine space. The matrix factorization constructed in {\itshape loc. cit.} is already $G$-equivariant and their formula (2.15) for the quantum dimensions agrees nicely with \refr{dim}.

\sssec{}{$K_{14}\sim Q_{10}$ in \cite{NRC}} In the notations of {\itshape loc. cit.} we take $G=\BZ/2\BZ$ acting by $(u,v,w)\mapsto (-u,v,-w-u^4)$ in case of version $1$ of $K_{14}$ or by $(u,v,w)\mapsto (u,-v,-w)$ in case of version $2$ of $K_{14}$ so that $f$ is of type $K_{14}$ and the singularity of $\hat{f}$ is of type $Q_{10}$. This was observed in \cite{BTW2} and once again a suitable affine space chart on $Y$ could be chosen. However, the matrix factorizations constructed in {\itshape loc. cit.} are not $G$-equivariant.

\sssec{}{$A_5\sim A_2\x A_2, Z_{13}\sim Q_{11}$ and $W_{13}\sim S_{11}$ in \cite{RW}} For $A_5\sim A_2\x A_2$ recall that we have just seen that $A_5\sim D_4$, so we have to check $D_4\sim A_2\x A_2$. Let $G=\BZ/3\BZ$ act in the notations of {\itshape loc. cit.} by $(y_1, y_2)\mapsto (\zeta_3 y_1, \zeta_3^{-1} y_2),$ where $\zeta_3$ is a nontrivial root of unity of degree $3$. In this case $f$ is of type $A_{2}\x A_2$ and the singularity of $\hat{f}$ is of type $D_{4}$. This was observed in \cite{BTW1}.

For $Z_{13}\sim Q_{11}$ one exposes the action easier for a different version of $Z_{13}$, namely $f=x^6_1x_2+x^3_2+x^2_3$, where $G=\BZ/2\BZ$ acts by $(x_1,x_2,x_3)\mapsto (-x_1,x_2,-x_3)$, so that the singularity of $\hat{f}$ is of type $Q_{11}$. This was observed in \cite{BTW2}.

For $W_{13}\sim S_{11}$ there is an action of $G=\BZ/2\BZ$ in the notations of {\itshape loc. cit.} given by $(y_1,y_2,y_3)\mapsto (-y_1,y_2,-y_3)$, so that $f$ is of type $W_{13}$ and the singularity of $\hat{f}$ is of type $S_{11}$. This was observed in \cite{BTW2}. The matrix factorization constructed in {\itshape loc. cit.} is already $G$-equivariant.

In all cases there are suitable affine space charts on $Y$.

\sssec{}{$E_{18}\sim Q_{12}$ and $E_{30}\sim Q_{18}$ in \cite{KRC}} In the notations of {\itshape loc. cit.} we take $G=\BZ/2\BZ$, which acts by $(x,y,z)\mapsto (-x,y,-z)$ for the first equivalence and by $(u,v,w)\mapsto (-u,v,-w+u^8)$ for the second. There are suitable affine space charts on $Y$ in both cases. 

\ssec{}{Nonexamples} $A_{11}\sim E_{6}, A_{17}\sim E_{7}$ and $A_{29}\sim E_{8}$ in \cite{CRCR}. The respective formulas (2.28), (2.41) and (2.51) for quantum dimensions in {\itshape loc. cit.} imply that the constructed matrix factorizations could not be obtained by our methods, even if we try to chain several equivalences of the form $\hat{f}\sim f$. Indeed by \refr{dim} we will always get a matrix factorization $X$ with $\dim_l(X)\dim_r(X)\in\BQ^{>0}$, while this is not the case in {\itshape loc. cit.} 

\ssec{}{Conjecture} Looking at these examples it is natural to pose the following conjecture:

\conj{conj}
1) Let $f\sim g$ be a (graded) orbifold equivalence of polynomials in either $2$ or $3$ variables given by a matrix factorization $X$, such that $\dim_l(X)\dim_r(X)\in\BQ^{>0}$. Then there exists a sequence of orbifold equivalences given by \reft{main}, such that $f\sim g$ is their compostion.

2) Let $f\sim g$ be a (graded) orbifold equivalence of singularities given by a matrix factorization $X$, such that $\dim_l(X)\dim_r(X)\in\BQ^{>0}$. Then there exists a sequence of orbifold equivalences, such that each equivalence is either of the type given by \reft{main} or is a Kn\"orrer periodicity (\cite[Proposition 1.2 (ii)]{CRCR}), so that $f\sim g$ is their compostion.
\econj



\begin{thebibliography}{9999}

\bibitem{BT} A.\,Basalaev, A.\,Takahashi, "Hochschild cohomology and orbifold Jacobian algebras associated to
invertible polynomials", J. Noncommut. Geom. 12;14(3):861--77, (2020).

\bibitem{BTW1} A.\,Basalaev, A.\,Takahashi, E.\,Werner, "Orbifold Jacobian algebras for invertible polynomials",
arXiv preprint: 1608.08962, (2016).

\bibitem{BTW2} A.\,Basalaev, A.\,Takahashi, E.\,Werner, "Orbifold Jacobian algebras for exceptional unimodal
singularities", Arnold Math. J., 3, 483--498, (2017).


\bibitem{BK} R.\,Bezrukavnikov, D.\,Kaledin, "McKay equivalence for symplectic resolutions of quotient singularities", Proc. Steklov Inst. Math. 246, no. 3 (2004), 13--33.

\bibitem{BGK} T.\,Bridgeland, A.\,King, M.\,Reid, "The McKay correspondence as
an equivalence of derived categories", J. Amer. Math. Soc. 14 (2001),
535--554.

\bibitem{BP} V.\,Baranovsky, J.\,Pecharich, "On equivalences of derived and singular categories", centr.eur.j.math. 8, (2010), 1--14.

\bibitem{CM} N.\,Carqueville, D.\,Murfet, "Adjunctions and defects in Landau-Ginzburg models", Adv. Math. 289 (2016), 480--566.

\bibitem{CRCR}  N.\,Carqueville, A.\,Ros Camacho, I.\,Runkel, "Orbifold equivalent
potentials", JPAA 220 2 (2016), 759--781.

\bibitem{CR} N.\,Carqueville, I.\,Runkel, "Rigidity and defect actions in Landau-Ginzburg models", Comm. Math. Phys. 310 (2012), 135--179.

\bibitem{CR2} N.\,Carqueville, I.\,Runkel, "Orbifold completion of defect bicategories", Quantum Topol. 7 (2016), no. 2, pp. 203--279

\bibitem{D} T.\,Dyckerhoff, "Compact generators in categories of matrix factorizations", Duke Math. J. 159
(2011), 223--274.

\bibitem{DM} T.\,Dyckerhoff, D.\,Murfet, "Pushing forward matrix factorisations", Duke Math. J. 162 (2013),
1249--1311.


\bibitem{H} M.\,Haiman, "Combinatorics, symmetric functions, and Hilbert schemes", Current developments
in mathematics, 2002, Int. Press, Somerville, MA (2003), 39--111.

\bibitem{KRC} T.\,Kluck, A.\,Ros Camacho, "Computational aspects of orbifold
equivalence", preprint, arXiv:1901.09019.

\bibitem{K} A.\,Kuznetsov, "Hyperplane sections and derived categories", Izv.: Math. 70:3, (2006), 447--547.

\bibitem{N}  T.\,Nevins, "Descent of coherent sheaves and complexes to geometric invariant theory quotients", Journal of
Algebra Volume 320, Issue 6, (2008), 2481--2495.

\bibitem{NRC} R.\,Newton, A.\,Ros Camacho, "Strangely dual orbifold equivalence
I", J. Sing. 14 (2016), 34--51.

\bibitem{O1} D.\,Orlov, "Triangulated categories of singularities and D-branes in Landau-Ginzburg models", Proc.
Steklov Inst. Math. 2004, no. 3(246), 227--248.

\bibitem{O2} D.\,Orlov, "Matrix factorizations for nonaffine LG-models", Math. Ann., 353, (2012), 95--108.

\bibitem{RW} A.\,Recknagel, P.\,Weinreb, "Orbifold equivalence: structure and new
examples", J. of Sing.
Vol. 17 (2018), 216--244.

\bibitem{R}  A.\,Ros Camacho, "On the Landau-Ginzburg/conformal field theory correspondence", Vertex Operator Algebras, Number Theory and Related Topics 753 (2020), Contemporary Mathematics AMS.

\bibitem{S} D.\,Shklyarov, "On Hochschild invariants of Landau-Ginzburg orbifolds", Adv. Theor. Math. Phys.,
24(1), 189--258, (2020).

\end{thebibliography}
\end{document}